\documentclass[11pt,reqno]{amsart}
\usepackage[T1]{fontenc}
\usepackage[english]{babel}
\usepackage{lmodern,microtype,mathtools,amssymb,amsthm}
\usepackage[mathcal]{euscript}
\usepackage[margin=1.1in]{geometry}
\usepackage{xcolor}
\usepackage[skins,breakable]{tcolorbox}
\definecolor{deepblue}{RGB}{30,78,140}
\definecolor{midblue}{RGB}{52,120,190}
\definecolor{urlpurple}{RGB}{120,81,169}
\definecolor{palepurple}{RGB}{236,230,246}
\colorlet{keyframe}{deepblue!80}
\colorlet{keyfill}{palepurple!70!white}
\usepackage{hyperref}
\hypersetup{colorlinks=true, linkcolor=deepblue, citecolor=midblue, urlcolor=urlpurple,
  pdftitle={A Stone-Weierstrass Theorem for Polish Spaces}, pdfauthor={Jose Iovino}}
\usepackage[capitalise]{cleveref}
\usepackage{enumitem}
\newtheoremstyle{rudin}{\topsep}{\topsep}{\itshape}{}{\bfseries}{.}{ }{\thmname{#1}\thmnumber{ #2}\thmnote{ (#3)}}
\newtheoremstyle{rudindef}{\topsep}{\topsep}{\upshape}{}{\bfseries}{.}{ }{\thmname{#1}\thmnumber{ #2}\thmnote{ (#3)}}
\theoremstyle{rudin}
\newtheorem{theorem}{Theorem}[section]
\newtheorem{keytheorem}[theorem]{Theorem}
\newtheorem{proposition}[theorem]{Proposition}
\newtheorem{lemma}[theorem]{Lemma}
\newtheorem{corollary}[theorem]{Corollary}
\newtheorem*{theoremA}{Theorem A}
\newtheorem*{theoremB}{Theorem B}
\theoremstyle{rudindef}
\newtheorem{definition}[theorem]{Definition}
\newtheorem{example}[theorem]{Example}
\newtheorem{remark}[theorem]{Remark}
\crefname{keytheorem}{Theorem}{Theorems}
\Crefname{keytheorem}{Theorem}{Theorems}
\tcolorboxenvironment{theoremA}{colback=keyfill, colframe=keyframe, boxrule=0.6pt, arc=6pt,
  left=6pt, right=6pt, top=2pt, bottom=2pt, before skip=10pt, after skip=10pt}
\tcolorboxenvironment{theoremB}{colback=keyfill, colframe=keyframe, boxrule=0.6pt, arc=6pt,
  left=6pt, right=6pt, top=2pt, bottom=2pt, before skip=10pt, after skip=10pt}
\tcolorboxenvironment{keytheorem}{colback=keyfill, colframe=keyframe, boxrule=0.4pt, arc=6pt,
  left=6pt, right=6pt, top=2pt, bottom=2pt, before skip=10pt, after skip=10pt}
\newcommand{\RR}{\mathbb R}
\newcommand{\QQ}{\mathbb Q}
\newcommand{\NN}{\mathbb N}
\newcommand{\MM}{\mathbb M}
\newcommand{\frag}{\mathcal F}
\newcommand{\cA}{\mathcal A}

\newcommand{\fm}{\mathfrak m}

\newcommand{\Loo}{L_{\omega_1\omega}}
\newcommand{\LR}{L_{\omega_1\omega}(\RR)}

\newcommand{\SFC}{S^{\frag}_C}
\newcommand{\HFC}{\mathcal H^{\frag}_C}
\newcommand{\open}[1]{\left\langle\!\left\langle #1\right\rangle\!\right\rangle}
\newcommand{\synt}[1]{#1}
\newcommand{\symb}[1]{\mathtt{#1}}
\newcommand{\unity}{\symb{1}}

\newcommand{\Sp}{\symb{p}}
\newcommand{\Sq}{\symb{q}}
\newcommand{\SQ}{\symb{Q}}
\newcommand{\Sneg}{\mathop{\synt{s\cdot}}}
\newcommand{\Splus}{\mathbin{\synt{+}}}
\newcommand{\Sor}{\mathop{\synt{\vee}}}
\newcommand{\Sand}{\mathop{\synt{\wedge}}}
\newcommand{\join}{\mathop{\synt{\bigvee}}}
\newcommand{\meet}{\mathop{\synt{\bigwedge}}}
\newcommand{\norm}[1]{\left\lVert#1\right\rVert_\infty}
\DeclareMathOperator{\BUC}{BUC}

\DeclareMathOperator{\Val}{Val}
\newcommand{\SE}{\symb{E}}
\title[Stone--Weierstrass for Polish spaces]{A Stone--Weierstrass Theorem for Polish Spaces,\\ Given by Logic}
\author{Jos\'e Iovino}
\address{Department of Mathematics, The University of Texas at San Antonio, San Antonio, TX 78249, USA}
\email{jose.iovino@utsa.edu}
\thanks{This version: \today. The most recent version of this paper is available at \url{https://utsa-quantum-math.github.io/resources/stone-weierstrass-polish-spaces/}; it may supersede the version deposited on arXiv.}
\date{November 6, 2026}
\begin{document}
\begin{abstract}
The Stone--Weierstrass theorem obtains every continuous function on a compact space as a uniform limit of functions built from a separating family by finitely many operations.
We extend the theorem to every separable metric space, and in particular to every Polish space, by admitting countable suprema and infima according to weights computed from the construction.
For a countable, uniformly bounded separating family whose supremum metric makes the space separable, the functions so built are exactly the bounded Lipschitz functions, and their uniform limits are exactly the bounded uniformly continuous functions.
The theorem was discovered, and first proved, in infinitary real-valued logic, where the weights are the a priori bounds of formulas and the uniform approximation is equivalent to an approximate form of Beth's definability theorem, which carries Reggio's correspondence from compact spaces to separable metric spaces.
The same bounds draw an exact line between the two forms of definability when the values of the coordinates in the models of the implicit definition form a set separable for the supremum metric: an implicitly defined predicate has an explicit definition precisely when the defining function is Lipschitz, and admits uniform approximation by sentences precisely when that function is uniformly continuous.
We also give a purely topological-analytic proof of the theorem.
\end{abstract}
\maketitle
\setcounter{tocdepth}{1}
\tableofcontents
\pagestyle{plain}

\section*{Introduction}

In this paper we extend the Stone--Weierstrass theorem from compact spaces to Polish spaces, and more generally to separable metric spaces.
The approximants are built from a countable separating family by rational linear combinations, maxima, minima, and countable suprema and infima, and a weight attached to each construction decides which countable suprema and infima are admitted; this single rule does the work that compactness does in the classical theorem.
The theorem was discovered, and first proved, in infinitary real-valued logic, whose a priori bounds are the weights and whose sentences define the approximants.
The proof by logic shows moreover that the uniform approximation is equivalent to an approximate form of Beth's definability theorem.
The statement, however, belongs to analysis, and for readers who do not work in model theory we give in the last section a purely topological-analytic proof, which requires no knowledge of logic.

Weierstrass proved in 1885 that every continuous function on a closed bounded interval is a uniform limit of polynomials~\cite{Weierstrass}, and half a century later Stone recognized that neither the interval nor the polynomials were essential~\cite{Stone1937,Stone1948}.
In its lattice form, Stone's theorem asserts that on a compact Hausdorff space every continuous real-valued function is a uniform limit of functions obtained from the rational constants and from an arbitrary separating family of continuous functions by finitely many sums, rational multiples, maxima and minima.
The theorem has, moreover, a logical counterpart: after Marra and Reggio had axiomatized the algebras of continuous functions from compact Hausdorff spaces to the unit interval, by finitary operations together with a single infinitary one~\cite{MarraReggio}, Reggio proved in 2021 that the Stone--Weierstrass theorem for these spaces is equivalent to the Beth definability property of the resulting infinitary equational logic~\cite{Reggio}.
Both results rest on compactness, and it is compactness that lets finitely many operations suffice.
Beyond compact spaces, finitely many operations do not suffice in general: on the countable discrete space, presented by the indicator functions of the points, every function built by finitely many operations is constant outside a finite set, and so stays at uniform distance at least $1$ from the function $n\mapsto(-1)^n$ (\cref{ex:discrete}).
A Stone--Weierstrass theorem for Polish spaces must therefore admit countable operations, and needs a principle that decides which of them to admit.

The extensions of the Stone--Weierstrass theorem beyond compact spaces have paid for the absence of compactness in one of three ways.
Some retain uniform approximation only in a weaker form.
Applied on compact subsets, the classical theorem yields approximation that is uniform on compact sets alone; applied to the one-point compactification of a locally compact space, it reaches only the functions that have a limit at infinity; and Buck's theorem concerns his strict topology, which is coarser than the topology of uniform convergence~\cite{Buck}, while Haydon showed that for the finer superstrict topology the theorem can hold only when that topology coincides with the strict one~\cite{Haydon}.
Others retain uniform approximation at the price of a stronger hypothesis on the family.
Hewitt, and Garrido and Montalvo after him, require the family to separate disjoint zero sets, or to satisfy conditions on the level sets of the function to be approximated~\cite{Hewitt,GarridoMontalvo1991,GarridoMontalvo1993}, and the passage through the Samuel compactification, which transfers the classical theorem to the bounded uniformly continuous functions on a metric space, requires the family to separate points of the compactification that the space does not contain~\cite{Samuel,Woods}.
The remaining extensions enlarge the operations themselves, but the enlargement is governed by an analytic property of the functions rather than by the way in which they are built.
Hager approximates every continuous function on a Lindel\"of space by a lattice-ordered group that separates points from closed sets, once the group is enlarged by local gluing, which admits a function according to how it behaves near each point~\cite{Hager1969,Hager2011}.
Mokhova, answering in this spirit a question of Cs\'asz\'ar~\cite{Csaszar}, obtains the bounded uniformly continuous functions on every uniform space from a vector lattice by envelopes followed by uniform limits~\cite{Mokhova,Raikov}; the envelopes, however, range over families of arbitrary cardinality, and a family is admitted only after its members have been shown to be equicontinuous.
None of these theorems obtains every bounded uniformly continuous function on a Polish space as a uniform limit of functions built from countably many separating functions by countable operations whose admission can be read off from the construction itself.
Theorem~A below achieves exactly this.

To state the theorem, fix a positive real number $R$, a nonempty set $X$, and functions $a_i\colon X\to[-R,R]$, for $i<\omega$, that separate the points of $X$.
These functions induce the \emph{coordinate metric}
\[
 d(x,y)=\sup_{i<\omega}|a_i(x)-a_i(y)|,
\]
the metric with respect to which Lipschitz continuity and uniform continuity are meant throughout.
The \emph{admissible functions} on $X$ are built together with \emph{weights}, by the following rules.
The constant function $1$ is admissible with weight $1$, and each function $a_i$ is admissible with weight $R$.
If $f$ and $g$ are admissible with weights $w$ and $w'$, and $s$ is a rational number, then $sf$, $f+g$, $\max(f,g)$ and $\min(f,g)$ are admissible with weights $|s|w$, $w+w'$, $\max(w,w')$ and $\max(w,w')$, respectively.
If $f_n$ is admissible with weight $w_n$ for every $n<\omega$, and $\sup_nw_n<\infty$, then $\sup_nf_n$ and $\inf_nf_n$ are admissible with weight $\sup_nw_n$.
The weight belongs to the construction rather than to the function, so the admission of each countable supremum or infimum is read off from the construction itself.

\begin{theoremA}
Let $R$ be a positive real number, let $X$ be a nonempty set, and let $(a_i)_{i<\omega}$ be functions from $X$ to $[-R,R]$ that separate the points of $X$, with coordinate metric $d$.
If the metric space $(X,d)$ is separable, then:
\begin{enumerate}[label=\textup{(\arabic*)}]
\item the admissible functions are exactly the bounded Lipschitz functions on $(X,d)$ (\cref{thm:representation,thm:an-representation});
\item the uniform closure of the set of admissible functions is the space $\BUC(X,d)$ of bounded uniformly continuous functions on $(X,d)$ (\cref{thm:approx,thm:an-approx}).
\end{enumerate}
\end{theoremA}

The operations are countable, and the only hypotheses on the family are the separation of points and the separability of the metric that the family induces.
Every separable metric space $(X,\rho)$ admits such a presentation, by the distances, truncated at $1$, to the points of a dense sequence, and the coordinate metric is then $\min(\rho,1)$, which has the same bounded uniformly continuous functions as $\rho$ (\cref{cor:polish}).
On every Polish space, therefore, every bounded uniformly continuous function for a compatible metric is a uniform limit of admissible functions built from the functions of the presentation; and when the space is compact, finitely many operations already suffice, by the lattice form of Stone's theorem (\cref{rem:compact}).

We write $\LR$ for infinitary real-valued logic, the $\Loo$ extension of Keisler's model theory for real-valued structures~\cite{Keisler,ADI}, in which the admission rule is a rule of syntax.
The \emph{coordinate vocabulary} has, for each $i<\omega$, a nullary predicate symbol $\Sp_i$ with a priori bound $R$, interpreted at a point $x$ of $X$ by the value $a_i(x)$, so that every sentence $\varphi$ for this vocabulary defines a function $x\mapsto\varphi^x$ on $X$.
The operations that build admissible functions are the connectives of $\LR$, and the weights are the a priori bounds that the formation rules attach to formulas; the admissible functions are therefore exactly the functions defined by sentences (\cref{lem:admissible}).

The second theorem carries Reggio's equivalence beyond compactness.
Let $\Sq$ be a new nullary predicate symbol, and let $B$ be a theory that implicitly defines $\Sq$ over the coordinate vocabulary, so that in every model of $B$ the value of $\Sq$ is a function $f_B$ of the values of the coordinate predicates.
The coordinate vocabulary has the \emph{approximate Beth property} when, for every such theory whose \emph{coordinate domain} $D_B$, the set of points of $[-R,R]^{\NN}$ that expand to models of $B$, is separable for the supremum metric and whose function $f_B$ is uniformly continuous on $D_B$, the value of $\Sq$ is a uniform limit, in all models of $B$, of the values of sentences for the coordinate vocabulary (\cref{def:abp}).

\begin{theoremB}
Let $R$ be a positive real number.
Then the coordinate vocabulary has the approximate Beth property (\cref{thm:abp}), and this property is equivalent to part~\textup{(2)} of Theorem~A for the number $R$ (\cref{prop:sw-beth,prop:beth-sw}).
\end{theoremB}

Both theorems rest on one observation about the a priori bound, the number that the formation rules attach to every formula and that bounds its values in every structure.
For the coordinate vocabulary, the a priori bound of a sentence, divided by the common bound $R$ of the coordinates, bounds the Lipschitz constant of the function that the sentence defines (\cref{lem:lip}).
The rule of syntax therefore keeps every admitted supremum among the Lipschitz functions, whereas a bound on values alone would admit the supremum of the functions $x\mapsto\max(-1,\min(1,nx))$ on $[-1,1]$, whose values are uniformly bounded although their a priori bounds are not (\cref{ex:overshoot}).
Conversely, every bounded Lipschitz function is a countable infimum of cones, each built from a countable supremum of coordinate differences: McShane's envelope~\cite{McShane} supplies the infimum, and the common a priori bound of the cones lets a single sentence of the original syntax define it (\cref{thm:representation}).
The sentences thus define the bounded Lipschitz functions, no more and no fewer.

For Theorem~B, the part played by Craig's interpolation theorem in the classical derivation of Beth's theorem~\cite{Beth,Craig} is taken by an approximate interpolation theorem for nullary vocabularies (\cref{thm:interpolation}).
The proof of the interpolation theorem translates the numerical entailment into infinitary propositional logic, through letters that record rational cuts with margins that make the translation robust, and separates the translated entailment by a positive formula obtained from a model existence game in the manner of Aubone and Iovino~\cite{AIgame}.
Threshold sentences, whose a priori bounds have a common bound, so that every countable operation on them is admitted, then translate the separator back.
Where the classical proof of interpolation invokes compactness, a position without a separator yields here, through a fair play of the game, a truth assignment that refutes the entailment.
Combined with a sentence that measures the distance to the graph of the implicitly defined function (\cref{lem:distance}), the interpolation theorem gives the approximate Beth property, and thereby part~(2) of Theorem~A.

The expressive power of $\LR$ reaches beyond the compactness theorem (\cref{ex:noncompact}), and the a priori bounds draw an exact line between the two forms of definability.
On a separable coordinate domain, an implicitly defined predicate has a definition by a single sentence if and only if the defining function is Lipschitz, and is uniformly approximable by sentences if and only if the function is uniformly continuous (\cref{lem:lip,thm:abp,cor:lipbeth}).
The square root of the absolute value of a coordinate realizes the distinction: the function is uniformly continuous but not Lipschitz, so the implicit definition admits arbitrarily accurate definitions by sentences but no exact one (\cref{ex:exactbeth}).
The a priori bound that makes Theorem~A exact is thus the same bound that makes Beth's theorem approximate; in the infinitary continuous logic of Ben Yaacov, Doucha, Nies and Tsankov, by contrast, whose countable connectives share a modulus of uniform continuity, numerical interpolation is exact~\cite[Corollary~6.4]{BDNT}.

For arbitrary countable vocabularies and countable theories, implicit definability yields a conclusion of another kind, which needs neither separability nor uniform continuity.
By the model existence theorem of Aubone, Due\~nez and Iovino~\cite{ADI}, the Henkin types of a theory, whose models are built from closed terms, form a Polish space.
When the theory implicitly defines a predicate, a Henkin type is determined by the values that it assigns to the atomic sentences in which that predicate does not occur, and the Lusin--Souslin theorem makes the set of these families of values Borel and recovers the Henkin type from them by a Borel map.
In particular, at every tuple of closed terms, the rational cuts of the implicitly defined predicate have Borel descriptions by atomic conditions in which that predicate does not occur (\cref{thm:borel}).

Once the theorem is known, analysis gives a second proof, presented in \cref{sec:analytic}.
Part~(1) follows from McShane's Lipschitz envelope, written as an infimum of countably many cones over a dense sequence, and part~(2) from the density of the Lipschitz functions among the bounded uniformly continuous functions, due to Georganopoulos for compact spaces and to Garrido and Jaramillo for metric spaces in general~\cite{Georganopoulos,GarridoJaramillo}.
The section uses only the admissible functions and their weights, and can be read immediately after this Introduction.

\subsection*{Organization of the paper}
In \cref{sec:framework} we recall the formulas and their a priori bounds.
Part~(1) of Theorem~A is proved in \cref{sec:representation}, the approximate interpolation theorem in \cref{sec:interpolation}, and Theorem~B, together with part~(2) of Theorem~A and the exact definitions that a Lipschitz hypothesis provides, in \cref{sec:equivalence}.
In \cref{sec:beth} we pass to arbitrary countable vocabularies and to the Borel descriptions of rational cuts.
In \cref{sec:analytic} we give the topological-analytic proof of Theorem~A, the presentation of separable metric spaces, and the examples that compare the theorem with the classical one.

\section{A rule of syntax bounds every formula}
\label{sec:framework}

In infinitary real-valued logic, the admissible functions of Theorem~A are the functions defined by sentences of the coordinate vocabulary (\cref{lem:admissible} below), and every estimate of the logical proof comes from the a priori bounds that the formation rules attach to formulas.
In this section we recall the syntax of the logic, the bounds and the conditions, following Aubone, Dueñez and Iovino~\cite{ADI}, whose structures are Keisler's general real-valued structures~\cite{Keisler}.
These structures interpret function symbols by arbitrary functions and predicate symbols by bounded real-valued functions, and the language has no equality symbol.
Every predicate symbol and every formula carries an individual bound, so that the logic has no fixed interval of truth values; the values $1$ and $-1$ are read as true and false.

A vocabulary $V=(\MM,F,P)$ prescribes nonempty sorts, function symbols with specified sorts, and predicate symbols $\Sp$ with individual a priori bounds $r_{\Sp}$.
The nullary logical unity symbol $\unity$ has value $1$ and a priori bound~$1$.
Every atomic formula $\Sp(\bar\tau)$ has a priori bound $r_{\Sp}$.
The first-order connectives are rational scalar multiplication, addition, binary maximum and supremum over a variable.
Their interpretations and a priori bounds are the following, for $s\in\QQ$:
\[
\begin{array}{c|c|c}
\text{formula}&\text{interpretation}&\text{a priori bound}\\ \hline
 \Sneg\varphi&s\varphi&|s|r_\varphi\\
 \varphi\Splus\psi&\varphi+\psi&r_\varphi+r_\psi\\
 \varphi\Sor\psi&\max(\varphi,\psi)&\max(r_\varphi,r_\psi)\\
 \synt{\sup_x}\varphi&\sup_x\varphi&r_\varphi
\end{array}
\]
The left column denotes syntax, and the middle column the corresponding operation on interpretations.
The minimum is the derived connective $-((-\varphi)\Sor(-\psi))$, with a priori bound $\max(r_\varphi,r_\psi)$.

The infinitary language $\Loo^V$ permits the formula $\join_n\varphi_n$ precisely when the formulas $\varphi_n$ have a common finite tuple of free variables and
\[
 \sup_n r_{\varphi_n}<\infty.
\]
The interpretation of $\join_n\varphi_n$ is the pointwise supremum, and the a priori bound is $\sup_n r_{\varphi_n}$.
The countable infimum $\meet_n\varphi_n$ is derived by negation and obeys the same bound rule.
The bounds are syntactic: they are fixed by the formation of the formula, not recalculated from interpretations on a particular class of structures.

A \emph{condition} is a strict rational inequality $\synt{\varphi>r}$ or $\synt{\varphi<r}$ on a sentence, and a \emph{theory} is a set of conditions.
Following Aubone, Due\~nez and Iovino~\cite{ADI}, the closed condition $\synt{\varphi\geq r}$ abbreviates the countable theory
\[
 \{\synt{\varphi>r-1/n}:n\geq1\},
\]
and the analogous convention applies to $\synt{\varphi\leq r}$.
The notation $\synt{\varphi=0}$ abbreviates the two closed conditions at zero; this equality concerns numerical values and adds no equality symbol to the vocabulary.

Real numerical constants can be represented while respecting the bounds.
Let a positive integer $M$ and a real number $c\in[-M,M]$ be given, and choose a sequence of rational numbers $t_n\in[-M,M]$ with supremum $c$.
Then
\[
 \underline c^{\,(M)}=\join_n(t_n\unity)
\]
is a sentence with value $c$ and a priori bound at most $M$.
The superscript records the chosen common bound and is omitted when that bound has been specified; no real scalar connective is added to the logic.

The expressive power of $\LR$ reaches beyond the compactness theorem, the tool behind the classical proof of Beth's theorem, already for nullary predicates.

\begin{example}[Beyond the compactness theorem]\label{ex:noncompact}
Let $(\Sp_n)_{n<\omega}$ be nullary predicate symbols with a priori bound~$1$.
The theory
\[
 \{\synt{\Sp_n\leq0}:n<\omega\}
 \ \cup\ \{\synt{\join_n\Sp_n>1/2}\}
\]
has no model, since the conditions $\synt{\Sp_n\leq0}$ force the value of $\join_n\Sp_n$ to be at most $0$.
Every finite subset of the theory, with the closed conditions expanded into strict ones, has a model: assign the value $1$ to a predicate whose index occurs in none of the finitely many conditions on individual predicates, and the value $0$ to all other predicates.
\end{example}

In the next section we determine what the countable operations admitted by the a priori bounds define.

\section{Sentences define exactly the bounded Lipschitz functions}
\label{sec:representation}

In this section we prove part~(1) of Theorem~A, the exact half of the Stone--Weierstrass theorem of this paper, in the language of sentences, which by \cref{lem:admissible} below define exactly the admissible functions.
On a set presented by countably many coordinate functions with a common bound, the a priori bound of a sentence, divided by the common bound, bounds the Lipschitz constant of the function that the sentence defines (\cref{lem:lip}), and, conversely, every bounded Lipschitz function is defined by a sentence when the coordinate metric is separable (\cref{thm:representation}).
The same estimate explains why the formation rule for countable suprema constrains a priori bounds rather than values (\cref{ex:overshoot}).

Let a positive real number $R$ be given.
The \emph{coordinate vocabulary} has one nonempty sort, no function symbols, and only the nullary predicate symbols $\unity$ and $\Sp_i$, where $i<\omega$ and $r_{\Sp_i}=R$; the sentences for this vocabulary are those of the infinitary language $\Loo$.
For the purpose of evaluating sentences, a structure for this vocabulary is determined by the point $z=(z_i)_{i<\omega}$ of $[-R,R]^\NN$ whose coordinates are the values of the predicates $\Sp_i$.
The sort does not affect the value of any sentence, since every quantifier over a variable is vacuous.
We write $\varphi(z)$ for the value of a sentence $\varphi$ at the point $z$, and we equip the cube $[-R,R]^\NN$ with the supremum metric
\[
 d_\infty(z,z')=\sup_i|z_i-z'_i|.
\]

\begin{lemma}\label{lem:lip}
Every sentence $\varphi$ for the coordinate vocabulary satisfies
\[
 |\varphi(z)-\varphi(z')|\leq \frac{r_\varphi}{R}\,d_\infty(z,z')
 \qquad(z,z'\in[-R,R]^\NN).
\]
\end{lemma}

\begin{proof}
The induction on formulas follows the arithmetic of the a priori bounds.
An atomic sentence $\Sp_i$ has Lipschitz constant $1=r_{\Sp_i}/R$, and the unity constant has Lipschitz constant~$0$.
Rational scalar multiplication multiplies a Lipschitz bound by the absolute value of the scalar, addition adds Lipschitz bounds, and the maximum of two functions preserves the larger bound.
For a countable join, the elementary inequality
\[
 |\sup_n u_n-\sup_n v_n|\leq\sup_n|u_n-v_n|
\]
gives the Lipschitz bound $R^{-1}\sup_n r_{\varphi_n}$, which is finite by the formation rule.
Quantifiers over variables do not change the function.
\end{proof}

A bound on the values of the terms of a countable supremum does not suffice for continuity.

\begin{example}[Bounded values do not suffice]\label{ex:overshoot}
Let $R=1$, let $X=[-1,1]$, let $a_0(x)=x$, and let $a_i(x)=0$ for $i\geq1$; the coordinate metric is $d(x,y)=|x-y|$.
For each positive integer $n$, the sentence
\[
 \chi_n=\synt{(-1)\cdot}\unity\Sor\bigl(\unity\Sand\synt{n\cdot}\Sp_0\bigr)
\]
defines the function $x\mapsto\max(-1,\min(1,nx))$, with values in $[-1,1]$, and has a priori bound $n$.
The supremum of these functions is $1$ for $x>0$ and $\max(-1,x)$ for $x\leq0$, a function that jumps from $0$ to $1$ at the origin.
By \cref{lem:lip}, no sentence defines this function, and indeed the formation rule excludes $\join_n\chi_n$, since the a priori bounds of the terms $\chi_n$ have no common bound.
\end{example}

We now interpret the coordinate predicates as a separating family on a set.
Let $X$ be a nonempty set, and let $a_i\colon X\to\RR$, for $i<\omega$, be functions that separate the points of $X$ and satisfy $|a_i(x)|\leq R$ for every $x\in X$ and every $i<\omega$.
At a point $x\in X$, the symbol $\Sp_i$ is interpreted by the value $a_i(x)$, and we write $\varphi^x$ for the value of a sentence $\varphi$ under this interpretation.
The \emph{coordinate metric} on $X$ is
\[
 d(x,y)=\sup_{i<\omega}|a_i(x)-a_i(y)|.
\]
Since the functions $a_i$ separate points, the function $d$ is a metric, bounded by $2R$.
No topology on $X$ is presupposed beyond the topology of this metric.

\begin{lemma}\label{lem:admissible}
Every sentence $\varphi$ for the coordinate vocabulary defines a function $x\mapsto\varphi^x$ that is admissible with weight $r_\varphi$, and every function that is admissible with weight $w$, in the sense of the Introduction, is defined by a sentence with a priori bound $w$.
\end{lemma}

\begin{proof}
Once the symbol $\unity$ is interpreted by the constant function $1$ and each symbol $\Sp_i$ by the function $a_i$, the formation rules of \cref{sec:framework}, with their a priori bounds, become the rules that build the admissible functions with their weights; the minimum and the countable infimum correspond to the derived connectives, which obey the same bound rules.
The first assertion follows by induction on sentences, since a supremum over a variable changes neither the function nor the a priori bound, and the second by induction on the construction of admissible functions.
\end{proof}

\begin{keytheorem}[Exact representation of bounded Lipschitz functions]\label{thm:representation}
Let $R$ be a positive real number, let $X$ be a nonempty set, and let $(a_i)_{i<\omega}$ be functions from $X$ to $[-R,R]$ that separate the points of $X$, with coordinate metric $d$.
If the metric space $(X,d)$ is separable, then the functions $x\mapsto\varphi^x$, for sentences $\varphi$ for the coordinate vocabulary, are exactly the bounded Lipschitz functions on $(X,d)$.
\end{keytheorem}

\begin{proof}
By \cref{lem:lip}, applied at the points $(a_i(x))_{i<\omega}$, every function defined by a sentence is bounded and Lipschitz on $(X,d)$.
For the converse, we write a bounded Lipschitz function as an infimum of cones based at a countable dense set.
Let $f\colon X\to\RR$ be bounded and Lipschitz.
Choose a positive integer $M>\norm f$, a positive rational number $L$ at least the Lipschitz constant of $f$, and a dense sequence $(x_n)_{n<\omega}$ in $(X,d)$.

Choose a positive integer $M_0\ge R$, and represent every constant $a_i(x_n)$ by a sentence with a priori bound at most $M_0$, as in \cref{sec:framework}.
The derived formula $|\varphi|=\varphi\Sor(-\varphi)$ has the same a priori bound as $\varphi$.
Hence
\[
 \delta_n=\join_{i<\omega}|\Sp_i-\underline{a_i(x_n)}|
\]
is a sentence with a priori bound at most $R+M_0$, and the value of $\delta_n$ at $x$ is $d(x,x_n)$.
Here $\delta_n$ is a name for a formula, not an additional connective.
Represent the constants $f(x_n)$ by sentences with a priori bound at most~$M$.
The formula
\[
 \varphi_f=\meet_{n<\omega}
           \bigl(\underline{f(x_n)}\Splus\synt{L\cdot}\delta_n\bigr)
\]
is permitted, since the inner formulas have the common a priori bound $M+L(R+M_0)$.

For every $x\in X$ and every $n<\omega$, the Lipschitz inequality gives
\[
 f(x)\leq f(x_n)+Ld(x,x_n).
\]
Conversely, for every positive number $\epsilon$, density gives an index $n$ with $d(x,x_n)<\epsilon$, and for that index we have
\[
 f(x_n)+Ld(x,x_n)\leq f(x)+2L\epsilon.
\]
Taking the infimum over $n$ and letting $\epsilon$ decrease to zero gives $\varphi_f^x=f(x)$.
\end{proof}

The analytic construction is the countable dense-set form of McShane's classical Lipschitz envelope~\cite{McShane}.
What the logic adds is an exact correspondence between syntax and analysis: the cones have a common a priori bound, so their infimum is again a sentence, and the sentences define the bounded Lipschitz functions and nothing else.
The logical proof of the approximation theorem also establishes its equivalence with approximate Beth definability, through the interpolation theorem proved next.

\section{A model existence game takes the place of compactness in interpolation}
\label{sec:interpolation}

In this section we prove approximate interpolation for countable compatible nullary vocabularies: a numerical entailment between sentences admits an interpolating sentence in the common vocabulary to every prescribed positive accuracy (\cref{thm:interpolation}), with no separability assumption on the space of valuations.
Classically, Beth's theorem is derived from interpolation, and the model-theoretic proof of interpolation rests on the compactness theorem; in infinitary real-valued logic, where the compactness theorem fails (\cref{ex:noncompact}), a model existence game takes the place of compactness, and the interpolation theorem becomes the logical engine of Theorem~B.
The proof translates the numerical problem into infinitary propositional logic, obtains a positive separator by a model existence game in the manner of Aubone and Iovino~\cite{AIgame}, and returns to the real-valued syntax through threshold sentences with a common a priori bound.

A \emph{nullary vocabulary} has one nonempty sort, no function symbols, and, besides the unity symbol $\unity$, only nullary predicate symbols, each with an a priori bound.
In a structure for a nullary vocabulary every quantifier over a variable is vacuous, so the value of a sentence depends only on the values of the predicate symbols.
We therefore identify a structure with a \emph{valuation}, an assignment $v$ of a real number $v(\Sp)$ to every predicate symbol $\Sp$ with $|v(\Sp)|\leq r_{\Sp}$ and $v(\unity)=1$, and we write $\varphi^v$ for the value of a sentence $\varphi$ at $v$.
Nullary vocabularies are \emph{compatible} if every predicate symbol common to them has the same a priori bound in each.
For compatible nullary vocabularies $V$ and $W$, a valuation for $V\cup W$ restricts to valuations for $V$, for $W$ and for $V\cap W$, and we evaluate a sentence for any of these vocabularies at a valuation for $V\cup W$ through the corresponding restriction.

\begin{keytheorem}[Approximate interpolation]\label{thm:interpolation}
Let $V$ and $W$ be countable compatible nullary vocabularies, let $\varphi$ be a sentence for $V$, and let $\psi$ be a sentence for $W$ such that $\varphi^v\leq\psi^v$ for every valuation $v$ for $V\cup W$.
Then for every real number $\epsilon>0$ there is a sentence $\theta$ for $V\cap W$ such that
\[
 \varphi^v-\epsilon\leq\theta^v\leq\psi^v
\]
for every valuation $v$ for $V\cup W$.
\end{keytheorem}

The proof occupies the rest of the section and has three ingredients: a perturbation lemma, a positive interpolation theorem for infinitary propositional logic, and a translation of sentences into propositional formulas.
A Boolean separator alone returns no real-valued sentence, and two devices bridge the gap.
Margins in the translation make the entailment robust, so that the values of the common predicates on the two sides are forced to be close (\cref{lem:robust}).
The threshold sentences substituted for the letters of the separator have a common a priori bound, so that every countable connective of the separator is admitted, and the positivity of the separator makes its numerical evaluation monotone, so that the resulting sentence is $1$ where $\varphi$ is large and $-1$ where $\psi$ is small (\cref{lem:saturation}).
The first ingredient is the induction of \cref{lem:lip}, stated for valuations.

\begin{lemma}\label{lem:perturbation}
Let $V$ be a nullary vocabulary, let $h\geq0$ be a real number, and let $v$ and $w$ be valuations for $V$ such that $|v(\Sp)-w(\Sp)|\leq h\,r_{\Sp}$ for every predicate symbol $\Sp$ of $V$.
Then every sentence $\varphi$ for $V$ satisfies $|\varphi^v-\varphi^w|\leq h\,r_\varphi$.
\end{lemma}

\begin{proof}
The induction on formulas in the proof of \cref{lem:lip} applies verbatim: the estimate holds for atomic sentences by hypothesis, and the connectives transform the bound $h\,r_\varphi$ exactly as they transform a priori bounds.
\end{proof}

The second ingredient is a positive interpolation theorem for infinitary propositional logic, which we prove by a model existence game.
A \emph{propositional formula} over a countable set of letters is built from the letters and the constants $\top$ and $\bot$ by negation and by countable conjunctions and disjunctions, and is evaluated at truth assignments to the letters in the usual way.
A propositional formula is \emph{positive} in a set $E$ of letters if, after negations have been pushed to the letters, no letter of $E$ occurs negated; a formula is \emph{positive} if no letter at all occurs negated after this transformation.

\begin{proposition}[Positive interpolation]\label{prop:lyndon}
Let $A$ and $B$ be propositional formulas that are positive in the set $E$ of letters common to $A$ and $B$, and suppose that every truth assignment that satisfies $A$ satisfies $B$.
Then there is a positive propositional formula $C$ whose letters belong to $E$ such that every truth assignment that satisfies $A$ satisfies $C$, and every truth assignment that satisfies $C$ satisfies $B$.
\end{proposition}

\begin{proof}
We play a model existence game in the manner of Aubone and Iovino~\cite{AIgame}, between Olivia, who challenges, and Josh, who answers.
Push the negations to the letters in $A$ and in $\neg B$, and write $B^-$ for the resulting form of $\neg B$; the letters of $E$ occur in $A$ only unnegated and in $B^-$ only negated.
A \emph{position} is a pair $(\Gamma,\Delta)$ of finite sets of formulas, the \emph{left} and the \emph{right commitments}, where $\Gamma$ consists of subformulas of $A$ and $\Delta$ of subformulas of $B^-$.
A \emph{separator} of the position is a positive propositional formula $C$ whose letters belong to $E$ such that every truth assignment that satisfies all members of $\Gamma$ satisfies $C$, and every truth assignment that satisfies all members of $\Delta$ falsifies $C$.
A position without a separator is \emph{inseparable}.
A separator of the position $(\{A\},\{B^-\})$ is the required formula, so it suffices to show that this position has a separator.

At each round Olivia names a commitment.
If the commitment is a conjunction $\bigwedge_n\varphi_n$, she also names an index $n$, and $\varphi_n$ is added to the same side; if the commitment is a disjunction $\bigvee_n\varphi_n$, Josh chooses an index $n$, and $\varphi_n$ is added to the same side.
Josh can keep the position inseparable.
Indeed, a separator of the position obtained by adding a conjunct separates the old position, since the conjunct follows from the conjunction.
If, for a disjunction $\bigvee_n\varphi_n$ among the left commitments, every choice of $n$ gave a position with a separator $C_n$, then $\bigvee_nC_n$ would separate the old position; for a disjunction among the right commitments, $\bigwedge_nC_n$ would.
Both formulas are positive, and their letters belong to $E$.
Moreover, an inseparable position contains no \emph{obstruction}: the constant $\bot$ on either side, a letter $\ell$ and $\neg\ell$ on the same side, or a letter $\ell\in E$ on the left together with $\neg\ell$ on the right.
Indeed, these obstructions are separated by $\bot$ or $\top$, according to the side, and by $\ell$, respectively.

Suppose that the position $(\{A\},\{B^-\})$ is inseparable.
There are countably many subformulas and countably many indices, so Olivia can schedule her challenges so that every conjunction that becomes a commitment is eventually named with every index, and every disjunction that becomes a commitment is eventually named.
Let Josh keep the position inseparable throughout, and let $\Gamma_\omega$ and $\Delta_\omega$ be the unions of the left and of the right commitments along the play.
Let $e_\omega$ be the truth assignment that makes a letter $p$ true exactly when either $p\in\Gamma_\omega$, or $p$ does not occur in $A$ and $\neg p\notin\Delta_\omega$.
Every literal in $\Gamma_\omega\cup\Delta_\omega$ is true at $e_\omega$.
For a literal on the left, the letter occurs in $A$, and the claim follows from the absence of obstructions on the left.
A literal on the right whose letter occurs in $A$ has a letter of $E$, which occurs in $B^-$ only negated, and the claim follows from the absence of the third kind of obstruction.
For a literal on the right whose letter does not occur in $A$, the claim follows from the definition of $e_\omega$ and the absence of obstructions on the right.
By induction on formulas, every member of $\Gamma_\omega\cup\Delta_\omega$ is true at $e_\omega$: the constant $\bot$ does not occur, every conjunct of a conjunction that occurs is added by the schedule, and some disjunct of a disjunction that occurs is added.
Hence $e_\omega$ satisfies $A$ and $B^-$, that is, $e_\omega$ satisfies $A$ and falsifies $B$, contrary to the hypothesis.
\end{proof}

The third ingredient translates sentences into propositional formulas.
Let $V$ be a countable nullary vocabulary.
For every predicate symbol $\Sp$ of $V$ with $r_{\Sp}>0$ and every rational number $u$ with $-r_{\Sp}<u<r_{\Sp}$, introduce a letter $\SE_{\Sp,u}$, the \emph{cut letters} of $V$.
A truth assignment $e$ to the cut letters \emph{decodes} to the valuation $\Val(e)$ given by
\[
 \Val(e)(\Sp)=\sup\bigl(\{-r_{\Sp}\}\cup\{u:\text{$e$ makes $\SE_{\Sp,u}$ true}\}\bigr)
\]
for every predicate symbol $\Sp$ with $r_{\Sp}>0$, by $\Val(e)(\Sp)=0$ when $r_{\Sp}=0$, and by $\Val(e)(\unity)=1$.
The supremum lies in $[-r_{\Sp},r_{\Sp}]$, so $\Val(e)$ is a valuation.
Every valuation $v$ arises in this way: the assignment that makes $\SE_{\Sp,u}$ true exactly when $v(\Sp)>u$ decodes to $v$, by density of the rational numbers.

\begin{lemma}\label{lem:cuts}
Let $V$ be a countable nullary vocabulary.
For every sentence $\varphi$ for $V$ and every rational number $t$ there are propositional formulas $\synt{[\varphi>t]}$ and $\synt{[\varphi<t]}$ over the cut letters of $V$ such that every truth assignment $e$ satisfies $\synt{[\varphi>t]}$ exactly when $\varphi^{\Val(e)}>t$, and satisfies $\synt{[\varphi<t]}$ exactly when $\varphi^{\Val(e)}<t$.
\end{lemma}

\begin{proof}
Let $\synt{[\varphi<t]}$ be the formula $\bigvee_{u\in\QQ,\,u<t}\neg\synt{[\varphi>u]}$, which has the required property once every $\synt{[\varphi>u]}$ has.
The formulas $\synt{[\varphi>t]}$ are defined by induction on $\varphi$.
For a predicate symbol $\Sp$ with $r_{\Sp}>0$, the formula $\synt{[\Sp>t]}$ is $\top$ if $t<-r_{\Sp}$, is $\bot$ if $t\geq r_{\Sp}$, and is otherwise the disjunction of the letters $\SE_{\Sp,u}$ with $u>t$; the definition of $\Val$ gives the required property.
For a symbol with a priori bound $0$ and for $\unity$, the formula is $\top$ or $\bot$ according to the constant value.
For rational scalars, $\synt{[s\cdot\varphi>t]}$ is $\synt{[\varphi>t/s]}$ if $s>0$, is $\synt{[\varphi<t/s]}$ if $s<0$, and is $\top$ or $\bot$ if $s=0$.
Finally,
\[
\begin{gathered}
 \synt{[\varphi\Splus\psi>t]}=\bigvee_{\substack{a,b\in\QQ\\ a+b>t}}\bigl(\synt{[\varphi>a]}\wedge\synt{[\psi>b]}\bigr),\\
 \synt{[\varphi\Sor\psi>t]}=\synt{[\varphi>t]}\vee\synt{[\psi>t]},
 \qquad
 \synt{[\join_n\varphi_n>t]}=\bigvee_n\synt{[\varphi_n>t]},
\end{gathered}
\]
and a supremum over a variable leaves the formula unchanged.
Each clause follows from the corresponding property of real numbers.
\end{proof}

We now fix the data of \cref{thm:interpolation}, with $U=V\cap W$.
For every predicate symbol $\Sp$ of $U$ with $r_{\Sp}>0$, choose a rational number $s_{\Sp}$ with $1\leq s_{\Sp}r_{\Sp}\leq2$.
Choose an integer $J>r_\varphi$, a rational number $\eta$ with $0<\eta\leq1$, and a finite set $T$ of rational numbers in $[-2,2]$ that contains $-2$ and $2$ and has gaps of length at most $\eta$, and let $g=3\eta J$.
Make two disjoint copies of the cut letters of $V\cup W$, the \emph{left} and the \emph{right} letters, and introduce new \emph{threshold letters} $\SE^+_{\Sp,t}$ and $\SE^-_{\Sp,t}$, for the predicate symbols $\Sp$ of $U$ with $r_{\Sp}>0$ and for $t\in T$.
We write $\synt{[\,\cdot\,]_{\mathrm L}}$ and $\synt{[\,\cdot\,]_{\mathrm R}}$ for the formulas of \cref{lem:cuts} written in the left and in the right letters, and, for a truth assignment $e$, we write $v_{\mathrm L}(e)$ and $v_{\mathrm R}(e)$ for the valuations decoded from the left and from the right letters.

For a rational number $a$, let $b=a-g$, and consider the propositional formulas
\[
\begin{aligned}
 A_a&=\synt{[\varphi>a]_{\mathrm L}}\wedge\bigwedge_{\Sp,t}\Bigl[\bigl(\neg\synt{[s_{\Sp}\cdot\Sp>t+\eta]_{\mathrm L}}\vee\SE^+_{\Sp,t}\bigr)\\
 &\qquad\wedge\bigl(\neg\synt{[s_{\Sp}\cdot\Sp<t-\eta]_{\mathrm L}}\vee\SE^-_{\Sp,t}\bigr)\Bigr],\\
 B_a&=\bigwedge_{n\geq1}\synt{[\psi>b-1/n]_{\mathrm R}}\vee\bigvee_{\Sp,t}\Bigl[\bigl(\SE^+_{\Sp,t}\wedge\neg\synt{[s_{\Sp}\cdot\Sp>t-\eta]_{\mathrm R}}\bigr)\\
 &\qquad\vee\bigl(\SE^-_{\Sp,t}\wedge\neg\synt{[s_{\Sp}\cdot\Sp<t+\eta]_{\mathrm R}}\bigr)\Bigr],
\end{aligned}
\]
where the formulas for $\varphi$ use only left letters of symbols of $V$, the formulas for $\psi$ use only right letters of symbols of $W$, and $\Sp$ and $t$ range over the threshold indices.
Thus $A_a$ says that $\varphi>a$ on the left and that every threshold relation covers the left values with margin $\eta$, and $B_a$ says that $\psi\geq b$ on the right or that some threshold relation is violated by the right values with margin $\eta$.
The letters common to $A_a$ and $B_a$ are the threshold letters, and both formulas are positive in them.

\begin{lemma}\label{lem:robust}
Let $a$ be a rational number.
Every truth assignment that satisfies $A_a$ satisfies $B_a$.
Consequently, there is a positive propositional formula $C_a$ in the threshold letters such that every truth assignment that satisfies $A_a$ satisfies $C_a$, and every truth assignment that satisfies $C_a$ satisfies $B_a$.
\end{lemma}

\begin{proof}
Let $e$ satisfy $A_a$, and suppose that no violation in $B_a$ occurs.
Let $\Sp$ be a predicate symbol of $U$ with $r_{\Sp}>0$, and write $P=v_{\mathrm L}(e)(\Sp)$ and $Q=v_{\mathrm R}(e)(\Sp)$.
If $s_{\Sp}P-s_{\Sp}Q>3\eta$, then the interval $[s_{\Sp}Q+\eta,\,s_{\Sp}P-\eta)$ has length greater than $\eta$ and lies in $[-2,2]$, so the interval contains a point $t$ of $T$; coverage gives $\SE^+_{\Sp,t}$, and the inequality $s_{\Sp}Q\leq t-\eta$ then gives a violation.
Symmetrically, if $s_{\Sp}Q-s_{\Sp}P>3\eta$, a point $t$ of $T$ in $(s_{\Sp}P+\eta,\,s_{\Sp}Q-\eta]$ gives $\SE^-_{\Sp,t}$ and a violation.
Hence $|P-Q|\leq3\eta/s_{\Sp}\leq3\eta\,r_{\Sp}$.

Let $w$ be the valuation for $V\cup W$ that agrees with $v_{\mathrm R}(e)$ on the symbols of $W$ and with $v_{\mathrm L}(e)$ on the symbols of $V$ outside $U$.
On the symbols of $V$, the valuations $w$ and $v_{\mathrm L}(e)$ differ by at most $3\eta\,r_{\Sp}$ at each symbol $\Sp$, so \cref{lem:perturbation} gives
\[
 \psi^{v_{\mathrm R}(e)}=\psi^w\geq\varphi^w\geq\varphi^{v_{\mathrm L}(e)}-3\eta\,r_\varphi>a-g=b.
\]
Thus $e$ satisfies every formula $\synt{[\psi>b-1/n]_{\mathrm R}}$, and hence $B_a$.
The formula $C_a$ is now given by \cref{prop:lyndon}.
\end{proof}

The positive formula $C_a$ is now converted into a sentence for $U$.
For a predicate symbol $\Sp$ of $U$ with $r_{\Sp}>0$ and for $t\in T$, the \emph{threshold sentences} are
\[
 \xi^+_{\Sp,t}=\synt{(-1)\cdot}\unity\Sor\Bigl(\unity\Sand\synt{\eta^{-1}\cdot}\bigl(s_{\Sp}\cdot\Sp\Splus\synt{(-t)\cdot}\unity\bigr)\Bigr),
 \qquad
 \xi^-_{\Sp,t}=\synt{(-1)\cdot}\unity\Sor\Bigl(\unity\Sand\synt{\eta^{-1}\cdot}\bigl(\synt{t\cdot}\unity\Splus\synt{(-s_{\Sp})\cdot}\Sp\bigr)\Bigr).
\]
Their values lie in $[-1,1]$, and their a priori bounds are at most $4/\eta$, since $s_{\Sp}r_{\Sp}\leq2$ and $|t|\leq2$.
Let $\theta_a$ be the sentence obtained from $C_a$, written without negations, by replacing each letter $\SE^\pm_{\Sp,t}$ by $\xi^\pm_{\Sp,t}$, the constants $\top$ and $\bot$ by $\unity$ and $\synt{(-1)\cdot}\unity$, and countable conjunctions and disjunctions by countable infima and suprema.
Each threshold sentence, and each sentence obtained at a node of $C_a$ after the substitution, has a priori bound at most $4/\eta$, however deep the formula $C_a$ is, so every countable infimum and supremum formed in this construction is permitted by the formation rules; the predicate symbols inside the threshold sentences keep their own a priori bounds.

The value of $\theta_a$ is computed by evaluating $C_a$ numerically: for a family $x$ of numbers in $[-1,1]$ indexed by the threshold letters, let $C_a(x)$ be the number obtained by interpreting each letter by the corresponding member of $x$, the constants by $1$ and $-1$, and conjunctions and disjunctions by infima and suprema.
Then $\theta_a^v=C_a(x)$, where $x$ is the family of values $(\xi^\pm_{\Sp,t})^v$.
Since $C_a$ is positive, $C_a(x)\leq C_a(y)$ whenever $x\leq y$ coordinatewise, and on families with values in $\{-1,1\}$ the number $C_a(x)$ is $1$ or $-1$ according as $C_a$ is true or false at the truth assignment that reads $1$ as true.

\begin{lemma}\label{lem:saturation}
Let $a$ be a rational number, and let $v$ be a valuation for $V\cup W$.
Then $-1\leq\theta_a^v\leq1$; moreover, $\theta_a^v=1$ if $\varphi^v>a$, and $\theta_a^v=-1$ if $\psi^v<a-g$.
\end{lemma}

\begin{proof}
The bounds hold because every operation in the construction of $\theta_a$ preserves $[-1,1]$.
Suppose that $\varphi^v>a$.
Let $e$ make the left letters describe $v$, as after the definition of $\Val$, and make $\SE^+_{\Sp,t}$ true exactly when $s_{\Sp}v(\Sp)>t+\eta$ and $\SE^-_{\Sp,t}$ true exactly when $s_{\Sp}v(\Sp)<t-\eta$.
Then $e$ satisfies $A_a$, hence $C_a$.
Whenever $e$ makes a threshold letter true, the corresponding threshold sentence has value $1$ at $v$; so the family of values of the threshold sentences dominates the $\{-1,1\}$-family of $e$, and monotonicity gives $\theta_a^v\geq1$.

Suppose instead that $\psi^v<a-g$.
Let $e$ make the right letters describe $v$, and make a threshold letter true exactly when the corresponding threshold sentence has value greater than $-1$ at $v$.
If $\SE^+_{\Sp,t}$ is true, then $s_{\Sp}v(\Sp)>t-\eta$, and if $\SE^-_{\Sp,t}$ is true, then $s_{\Sp}v(\Sp)<t+\eta$; so no violation occurs, and $e$ falsifies $B_a$ because $\psi^v<b$.
Hence $e$ falsifies $C_a$.
The family of values of the threshold sentences is dominated by the $\{-1,1\}$-family of $e$, and monotonicity gives $\theta_a^v\leq-1$.
\end{proof}

\begin{proof}[Proof of \cref{thm:interpolation}]
Let $\epsilon>0$ be given.
Choose a positive integer $M$ greater than $r_\varphi$ and $r_\psi$, the integer $J$ and the rational number $\eta$ of the construction above, with $\eta$ so small that $g=3\eta J<\epsilon/2$, and a positive integer $N$ such that $\delta=2M/N$ satisfies $\delta\leq\epsilon/2$.
For $1\leq j\leq N$ let $c_j=-M+j\delta$, and let
\[
 \theta=\synt{(-M-g)\cdot}\unity\Splus\synt{(\delta/2)\cdot}\bigl((\theta_{c_1}\Splus\unity)\Splus\cdots\Splus(\theta_{c_N}\Splus\unity)\bigr),
\]
a sentence for $U$.
Fix a valuation $v$ for $V\cup W$, and let $w_j=(\theta_{c_j}^v+1)/2\in[0,1]$, so that $\theta^v=-M-g+\delta\sum_jw_j$.
By \cref{lem:saturation}, $w_j=1$ whenever $c_j<\varphi^v$; there are at least $(\varphi^v+M)/\delta-1$ such indices, so $\theta^v\geq\varphi^v-\delta-g\geq\varphi^v-\epsilon$.
By \cref{lem:saturation} again, $w_j=0$ unless $c_j\leq\psi^v+g$; there are at most $(\psi^v+g+M)/\delta$ such indices, so $\theta^v\leq\psi^v$.
\end{proof}

Interpolation is the engine of Beth's theorem, and in the next section the approximate interpolation theorem yields the approximate Beth property and thereby part~(2) of Theorem~A.

\section{Approximate Beth definability is equivalent to the new Stone--Weierstrass theorem}
\label{sec:equivalence}

In this section we prove part~(2) of Theorem~A by logic.
Stated for sentences, as \cref{thm:approx} below, part~(2) is equivalent to an approximate form of Beth's definability theorem for the coordinate vocabulary, and we prove the approximate form from the interpolation theorem of \cref{sec:interpolation}.
The equivalence is the analogue, for spaces that need not be compact, of Reggio's equivalence between the Stone--Weierstrass theorem and the Beth definability property of an infinitary equational logic~\cite{Reggio}.
Beth's theorem holds in the logic in exact form under a Lipschitz hypothesis (\cref{cor:lipbeth}), and in approximate form for uniformly continuous definitions (\cref{thm:abp}), which need not admit exact ones (\cref{ex:exactbeth}).

Under the hypotheses of \cref{thm:representation}, we write $\BUC(X,d)$ for the space of bounded uniformly continuous real-valued functions on $(X,d)$ and $\mathcal E(X)$ for the set of functions $x\mapsto\varphi^x$ defined by sentences for the coordinate vocabulary, which by \cref{lem:admissible} is the set of admissible functions.
Uniform closures are taken in the space of all bounded real-valued functions on $X$.
Part~(2) of Theorem~A then reads as follows; the proof is given at the end of the section.

\begin{keytheorem}[Stone--Weierstrass theorem for separable metric spaces]\label{thm:approx}
Let $R$ be a positive real number, let $X$ be a nonempty set, and let $(a_i)_{i<\omega}$ be functions from $X$ to $[-R,R]$ that separate the points of $X$, with coordinate metric $d$.
If the metric space $(X,d)$ is separable, then the uniform closure of $\mathcal E(X)$ is $\BUC(X,d)$.
\end{keytheorem}

To formulate the approximate form of Beth's theorem, let a positive real number $R$ be given, and recall the coordinate vocabulary of \cref{sec:representation}.
Let $\Sq$ be a nullary predicate symbol, not in the coordinate vocabulary, with a priori bound $r_{\Sq}>0$.
A structure for the coordinate vocabulary expanded by $\Sq$ is determined by a point $z$ of $[-R,R]^\NN$ and a value $t\in[-r_{\Sq},r_{\Sq}]$ of $\Sq$, and we write $(z,t)$ for the structure.
For a theory $B$ for the expanded vocabulary, let the \emph{coordinate domain} $D_B$ of $B$ be the set of points $z$ such that $(z,t)$ is a model of $B$ for some $t$.
The theory $B$ \emph{implicitly defines} $\Sq$ if $t=t'$ whenever $(z,t)$ and $(z,t')$ are models of $B$; in that case the \emph{implicitly defined function} $f_B\colon D_B\to[-r_{\Sq},r_{\Sq}]$ assigns to each point $z$ of $D_B$ the unique such value $t$.
We equip $D_B$ with the supremum metric $d_\infty$.

\begin{definition}\label{def:abp}
The coordinate vocabulary has the \emph{approximate Beth property} if the following holds for every nullary predicate symbol $\Sq$ as above and every theory $B$ that implicitly defines $\Sq$ such that the metric space $(D_B,d_\infty)$ is separable and the function $f_B$ is uniformly continuous on that space: for every real number $\epsilon>0$ there is a sentence $\theta$ for the coordinate vocabulary such that $|t-\theta(z)|\leq\epsilon$ for every model $(z,t)$ of $B$.
\end{definition}

The hypothesis of uniform continuity is necessary: if such sentences exist for every $\epsilon$, then $f_B$ is a uniform limit of the functions $z\mapsto\theta(z)$, which are Lipschitz by \cref{lem:lip}, so $f_B$ is uniformly continuous.

The passage from the approximate Beth property back to \cref{thm:approx}, and the proof of that property from interpolation, use sentences that measure the distance to a graph.
On $[-R,R]^\NN\times[-r_{\Sq},r_{\Sq}]$ we use the metric $\rho((z,t),(z',u))=\max(d_\infty(z,z'),|t-u|)$.

\begin{lemma}\label{lem:distance}
Let $G$ be a nonempty separable subset of $([-R,R]^\NN\times[-r_{\Sq},r_{\Sq}],\rho)$.
Then there is a sentence $\gamma_G$ for the coordinate vocabulary expanded by $\Sq$ whose value at every structure $(z,t)$ is the distance $\rho((z,t),G)$.
\end{lemma}

\begin{proof}
Choose a dense sequence $((y_n,u_n))_{n<\omega}$ in $G$ and a positive integer $M_0\geq\max(R,r_{\Sq})$, and represent the constants $y_{n,i}$ and $u_n$ by sentences with a priori bound at most $M_0$, as in \cref{sec:framework}.
The sentence
\[
 \gamma_G=\meet_{n<\omega}\Bigl(\join_{i<\omega}|\Sp_i-\underline{y_{n,i}}|\ \Sor\ |\Sq-\underline{u_n}|\Bigr)
\]
is permitted, since the inner formulas have the common a priori bound $M_0+\max(R,r_{\Sq})$.
The value of $\gamma_G$ at $(z,t)$ is $\inf_n\rho((z,t),(y_n,u_n))$, which equals $\rho((z,t),G)$ by density.
\end{proof}

\begin{proposition}\label{prop:sw-beth}
Let $R$ be a positive real number.
If the conclusion of \cref{thm:approx} holds for $R$, then the coordinate vocabulary for $R$ has the approximate Beth property.
\end{proposition}

\begin{proof}
Let $\Sq$ and $B$ be as in \cref{def:abp}, and let $\epsilon>0$ be given.
If $D_B$ is empty, then $B$ has no model and every sentence will do.
Otherwise, apply the conclusion of \cref{thm:approx}, which holds for $R$ by hypothesis, to the set $D_B$ and the coordinate functions $z\mapsto z_i$, which take values in $[-R,R]$, separate points and have coordinate metric $d_\infty$.
The function $f_B$ is bounded by $r_{\Sq}$ and uniformly continuous, so there is a sentence $\theta$ with $|f_B(z)-\theta(z)|\leq\epsilon$ for every $z\in D_B$.
For every model $(z,t)$ of $B$ we have $t=f_B(z)$, which gives the conclusion.
\end{proof}

\begin{proposition}\label{prop:beth-sw}
Let $R$ be a positive real number, and suppose that the coordinate vocabulary for $R$ has the approximate Beth property.
Let $X$ be a nonempty set, and let $(a_i)_{i<\omega}$ be functions from $X$ to $[-R,R]$ that separate the points of $X$, with coordinate metric $d$ such that $(X,d)$ is separable.
Then every function in $\BUC(X,d)$ is a uniform limit of functions $x\mapsto\varphi^x$ defined by sentences $\varphi$ for the coordinate vocabulary.
\end{proposition}

\begin{proof}
Let $f\in\BUC(X,d)$ and $\epsilon>0$ be given, and choose a positive integer $r_{\Sq}>\norm f$ as the a priori bound of $\Sq$.
The map $\iota\colon X\to[-R,R]^\NN$ given by $\iota(x)=(a_i(x))_{i<\omega}$ is an isometric embedding of $(X,d)$ into $([-R,R]^\NN,d_\infty)$.
The cube is complete for $d_\infty$, so the uniformly continuous function $f\circ\iota^{-1}$ on $\iota(X)$ extends to a uniformly continuous function $F$ on the closure $K$ of $\iota(X)$, with values in $[-\norm f,\norm f]$.
The graph $G$ of $F$ is closed and separable, and \cref{lem:distance} gives the sentence $\gamma_G$.
Let $B$ be the theory $\{\synt{\gamma_G\leq0}\}$.
The models of $B$ are the structures $(z,t)$ with $(z,t)\in G$, so $B$ implicitly defines $\Sq$, with $D_B=K$ and $f_B=F$.
The space $K$ is separable and $F$ is uniformly continuous, so the approximate Beth property gives a sentence $\theta$ with $|F(z)-\theta(z)|\leq\epsilon$ for every $z\in K$.
At $z=\iota(x)$ this inequality reads $|f(x)-\theta^x|\leq\epsilon$.
\end{proof}

Together with \cref{lem:lip}, \cref{prop:beth-sw} yields \cref{thm:approx}.
We now prove the approximate Beth property without \cref{thm:approx}, from interpolation.

\begin{keytheorem}[Approximate Beth definability]\label{thm:abp}
For every positive real number $R$, the coordinate vocabulary has the approximate Beth property.
\end{keytheorem}

\begin{proof}
Let $\Sq$ and $B$ be as in \cref{def:abp}, and let $\epsilon>0$ be given.
If $D_B$ is empty, then $B$ has no model and every sentence will do.
Otherwise, as in the proof of \cref{prop:beth-sw}, the function $f_B$ extends to a uniformly continuous function $F$ on the closure $K$ of $D_B$, and the graph $G$ of $F$ is nonempty and separable.
Let $\Sq'$ be a new nullary predicate symbol with a priori bound $r_{\Sq}$, and let $\gamma'_G$ be the sentence obtained from $\gamma_G$ by replacing $\Sq$ with $\Sq'$.
Choose a real number $\delta_0>0$ such that $|F(y)-F(y')|<\epsilon/2$ whenever $y,y'\in K$ and $d_\infty(y,y')<\delta_0$, a positive integer $N\geq2r_{\Sq}/\delta_0$, and a rational number $\lambda$ with $\epsilon/2\leq\lambda\leq3\epsilon/4$.
Consider the sentences
\[
 \varphi=\Sq\Splus\synt{(-N)\cdot}\gamma_G,
 \qquad
 \psi=\Sq'\Splus\synt{N\cdot}\gamma'_G\Splus\synt{\lambda\cdot}\unity
\]
for the compatible nullary vocabularies obtained by adding $\Sq$ and $\Sq'$, respectively, to the coordinate vocabulary.

We claim that $\varphi\leq\psi$ at every valuation, that is, at every point $z$ of the cube and all values $t,t'\in[-r_{\Sq},r_{\Sq}]$ of $\Sq$ and $\Sq'$.
Let $y,y'\in K$ be given, and let $\alpha=\rho((z,t),(y,F(y)))$ and $\alpha'=\rho((z,t'),(y',F(y')))$, so that $d_\infty(y,y')\leq\alpha+\alpha'$.
If $\alpha+\alpha'<\delta_0$, then, since $N\geq1$,
\[
 t-t'\leq|t-F(y)|+|F(y)-F(y')|+|F(y')-t'|\leq\alpha+\alpha'+\epsilon/2\leq N(\alpha+\alpha')+\lambda.
\]
If $\alpha+\alpha'\geq\delta_0$, then $t-t'\leq2r_{\Sq}\leq N\delta_0\leq N(\alpha+\alpha')$.
Hence $t-N\alpha\leq t'+N\alpha'+\lambda$ for all $y,y'\in K$, and taking the infimum of $\alpha$ over $y$ and of $\alpha'$ over $y'$ gives $t-N\gamma_G(z,t)\leq t'+N\gamma'_G(z,t')+\lambda$, which is the claim.

Applied with $\epsilon/4$, \cref{thm:interpolation} gives a sentence $\theta$ for the coordinate vocabulary with $\varphi-\epsilon/4\leq\theta\leq\psi$ at every valuation.
Let $(z,t)$ be a model of $B$.
Then $z\in D_B$, $t=F(z)$ and $\gamma_G(z,t)=0$, and evaluating at the valuation in which $\Sq'$ also takes the value $t$ gives
\[
 t-\epsilon/4\leq\theta(z)\leq t+\lambda\leq t+3\epsilon/4.
\]
Hence $|t-\theta(z)|\leq\epsilon$.
\end{proof}

The approximate Beth property now yields part~(2) of Theorem~A.

\begin{proof}[Proof of \cref{thm:approx}]
By \cref{lem:lip}, every function in $\mathcal E(X)$ is Lipschitz, so every uniform limit of functions in $\mathcal E(X)$ is bounded and uniformly continuous.
Conversely, the coordinate vocabulary for $R$ has the approximate Beth property by \cref{thm:abp}, and \cref{prop:beth-sw} shows that every function in $\BUC(X,d)$ is a uniform limit of functions in $\mathcal E(X)$.
\end{proof}

An implicit definition by a uniformly continuous function, approximable by \cref{thm:abp}, need not be explicit.

\begin{example}[An implicit definition beyond Lipschitz functions]\label{ex:exactbeth}
Let $R=1$ and $r_{\Sq}=1$, let $K$ be the set of points $(x,0,0,\dots)$ of the cube with $x\in[-1,1]$, and let $F(x,0,0,\dots)=\sqrt{|x|}$.
The set $K$ is separable, and the graph $G$ of $F$ is closed.
The theory $B=\{\synt{\gamma_G\leq0}\}$ implicitly defines $\Sq$, with $D_B=K$ and $f_B=F$, and $F$ is uniformly continuous; so \cref{thm:abp} applies.
No sentence $\theta$ for the coordinate vocabulary satisfies $\theta(z)=t$ for every model $(z,t)$ of $B$: such a sentence would satisfy $\theta(x,0,0,\dots)=\sqrt{|x|}$, which is not Lipschitz, contrary to \cref{lem:lip}.
\end{example}

Under a Lipschitz hypothesis, by contrast, \cref{thm:representation} gives an exact definition by a single sentence, and together with \cref{lem:lip} the conclusion characterizes exact definability on separable coordinate domains: an implicitly defined predicate is explicitly definable precisely when the defining function is Lipschitz.

\begin{corollary}[Lipschitz Beth definability]\label{cor:lipbeth}
Let $B$ be a theory for the coordinate vocabulary expanded by $\Sq$ that implicitly defines $\Sq$.
If the metric space $(D_B,d_\infty)$ is separable and the function $f_B$ is Lipschitz on that space, then there is a sentence $\varphi$ for the coordinate vocabulary such that every model of $B$ satisfies $\synt{\Sq-\varphi=0}$.
\end{corollary}

\begin{proof}
If $D_B$ is empty, the theory $B$ has no model, and every sentence will do.
Otherwise, apply \cref{thm:representation} to the set $D_B$ and the coordinate functions $z\mapsto z_i$.
The resulting sentence takes the value $f_B(z)$, and hence the value of $\Sq$, in every model $(z,t)$ of $B$.
\end{proof}

The sentence of \cref{cor:lipbeth} belongs to the original real-valued syntax, with rational scalars only and with a common bound on the a priori bounds at every countable join and infimum.
Approximation and definability are thus bound together for the coordinate vocabulary; for arbitrary countable vocabularies, implicit definitions find in the next section a description of another kind.

\section{Beyond the coordinate vocabulary, implicit definitions have Borel descriptions}
\label{sec:beth}

In this section we leave the coordinate vocabulary for arbitrary countable vocabularies, and the question changes accordingly: implicit definability is no longer accompanied by the separability and uniform continuity of \cref{sec:equivalence}, and the conclusion is a Borel description rather than an approximation by sentences.
Descriptive set theory takes the place of the compactness theorem.
On a compact space of models, a continuous injective reduct map would have a continuous inverse on its image; on the Polish spaces of Henkin types, the Lusin--Souslin theorem supplies a Borel inverse, and the rational cuts of an implicitly defined predicate therefore have Borel descriptions.
We begin with the model existence theorem of Aubone, Due\~nez and Iovino~\cite{ADI}.

Let $V$ be a countable vocabulary, let $C$ be a set of new constants containing countably many constants of every sort, and let $\frag\subseteq\Loo^V$ be a countable fragment.
Let $V_C$ be the vocabulary obtained by adding the constants of $C$ to $V$; the fragment generated by $\frag$ over $V_C$ has a countable set $\SFC$ of sentences.
Formal types are points of the compact metrizable space
\[
 \prod\SFC=\prod_{\varphi\in\SFC}[-r_\varphi,r_\varphi],
\]
and for a formal type $\fm$ we write $\varphi^\fm$ for the coordinate of $\fm$ at $\varphi$.
Henkin types are the formal types that satisfy the recursive equations for the connectives and the witness equations
\[
 (\synt{\sup_x}\varphi)^\fm
   =\sup_{\tau}\varphi(\tau)^\fm,
 \qquad
 (\join_n\varphi_n)^\fm=\sup_n\varphi_n^\fm,
\]
where $\tau$ ranges over the closed terms of the required sort.
The set of Henkin types is denoted by $\HFC$, and for a theory $B$ the set $\open B$ consists of the formal types that satisfy every condition of~$B$.

\begin{proposition}[Henkin model existence]\label{prop:existence}
Let $V$, $C$ and $\frag$ be as above.
Every member of $\HFC$ is realized by a term structure whose elements in each sort are the closed terms of that sort.
The set $\HFC$ is a $G_\delta$ subset of $\prod\SFC$, and for every theory $B$ on $\SFC$ the set $\HFC\cap\open B$ is \v{C}ech-complete and Polish.
\end{proposition}

\begin{proof}
The first two assertions are the Henkin realization lemma and the topology theorem of Aubone, Due\~nez and Iovino~\cite{ADI}.
For realization, interpret the function symbols by term formation and each predicate symbol on a tuple of closed terms by the prescribed atomic value.
Induction on formulas, with the Henkin witness equations at the two supremum steps, gives all the prescribed values.
No quotient is required, because the language has no equality symbol.

For the topological assertion, the Henkin equation for a countable supremum $\sigma=\join_n\varphi_n$ is the conjunction of the closed inequalities $\varphi_n\le\sigma$ and the requirement
\[
 \bigcap_{k\geq1}\bigcup_n
 \{\fm:\sigma^\fm<\varphi_n^\fm+1/k\},
\]
which is a $G_\delta$ condition.
The equation for a supremum over a variable has the same form, indexed by closed terms, and the other Henkin equations are closed.
There are countably many equations, so the set $\HFC$ is $G_\delta$.
Each strict condition defines an open set, and there are only countably many conditions on $\SFC$, so the set $\HFC\cap\open B$ is $G_\delta$ in the compact metrizable space $\prod\SFC$, and is therefore \v{C}ech-complete and Polish.
\end{proof}

To formulate uniqueness, let $V'=V\cup\{\SQ\}$, where $\SQ$ is a new predicate symbol with prescribed sorts and finite a priori bound; the vocabulary $V'$ has the same sorts and function symbols as $V$.
A $V'$-theory $B$ \emph{implicitly defines} $\SQ$ over $V$ if any two models of $B$ with the same $V$-reduct have the same interpretation of $\SQ$, as a function on the shared sorts.
A \emph{Borel cut description} is a well-founded expression formed from strict rational conditions on atomic $V_C$-sentences by countable unions, countable intersections and complements.
The operations refer to sets of models, and a Borel cut description is not asserted to be a formula of $\Loo^{V_C}$.
Fix a countable fragment that contains the sentences of $B$.
A \emph{named Henkin model} of $B$ is the term structure, given by \cref{prop:existence}, of a Henkin type over $V'_C$ that satisfies $B$, and the \emph{reduct atomic type} of such a model is the family of the values of the atomic $V_C$-sentences in the model.

\begin{keytheorem}[Beth definability of rational cuts]\label{thm:borel}
Let $V$ be a countable vocabulary, and let $B$ be a countable $V'$-theory that implicitly defines $\SQ$ over $V$.
For every tuple $\bar\tau$ of closed $V_C$-terms of the required sorts and every rational number $r$, the condition $\synt{\SQ(\bar\tau)>r}$ has, on the named Henkin models of $B$, a Borel cut description that uses only atomic $V_C$-conditions.
The set of reduct atomic types of these models is Borel as well.
\end{keytheorem}

\begin{proof}
We use Henkin realization to prove that the reduct map is injective, and then apply the Lusin--Souslin theorem.
Form the Henkin space over $V'_C$ for the fixed fragment.
Let $\mathcal P$ be the Polish space of Henkin types that satisfy $B$, provided by \cref{prop:existence}.
Let $\cA_V$ be the countable set of atomic $V_C$-sentences, and let $Y=\prod_{\alpha\in\cA_V}[-r_\alpha,r_\alpha]$.
The coordinate restriction
\[
 \pi\colon\mathcal P\longrightarrow Y
\]
is continuous.

Suppose that two members $\fm$ and $\fm'$ of $\mathcal P$ have the same image under $\pi$.
The term models of $\fm$ and $\fm'$ have the same closed-term universes and the same interpretations of the function symbols.
Agreement on the atomic $V_C$-sentences gives the same reduct predicates, so the two term models have the same $V$-reduct.
Since $B$ implicitly defines $\SQ$, the two term models interpret $\SQ$ in the same way.
They therefore agree on every atomic $V'_C$-sentence, and the recursive Henkin equations give $\fm=\fm'$.
Hence the map $\pi$ is injective.

The Lusin--Souslin theorem states that the image of a Borel subset of a Polish space under an injective continuous map into a Polish space is Borel~\cite[Theorem~15.1]{Kechris}.
Applied to $\mathcal P$, the theorem shows that the set $\pi(\mathcal P)$ is Borel in~$Y$.
Applied to the open subset
\[
 \{\fm\in\mathcal P:\SQ(\bar\tau)^\fm>r\}
\]
of $\mathcal P$, the theorem shows that the image of this subset is Borel in~$Y$.
The topology of $Y$ has a countable subbasis of cylinders given by strict rational atomic inequalities, so every Borel subset of $Y$ has a description by countable Boolean operations on these cylinders.
Pulling the resulting description back along $\pi$ proves the required equivalence on $\mathcal P$.
\end{proof}

The descriptions of \cref{thm:borel} concern every named closed-term instance of the predicate, with the same constants used to code the reduct.
No elimination of these naming constants, and no definition by a single real-valued formula, is claimed.
Henkin realizability supplies actual models, the Lusin--Souslin theorem supplies Borel inversion, and neither step uses the compactness theorem.

\section{Theorem A has a topological-analytic proof}
\label{sec:analytic}

The logical proof of part~(2) of Theorem~A in \cref{sec:equivalence} passes through interpolation and the approximate Beth property.
Once the theorem is known, analysis alone proves both parts: part~(1) follows from McShane's Lipschitz envelope~\cite{McShane}, taken over a dense sequence, and part~(2) from the density of the Lipschitz functions among the bounded uniformly continuous functions~\cite{Georganopoulos,GarridoJaramillo}.
The section uses only the admissible functions and their weights, defined in the Introduction, and is independent of \cref{sec:framework,sec:representation,sec:interpolation,sec:equivalence,sec:beth}.

\begin{lemma}\label{lem:an-weights}
Let $R$ be a positive real number, let $X$ be a nonempty set, and let $(a_i)_{i<\omega}$ be functions from $X$ to $[-R,R]$ that separate the points of $X$, with coordinate metric $d$.
If $f$ is admissible with weight $w$, then $|f(x)|\leq w$ and $|f(x)-f(y)|\leq(w/R)\,d(x,y)$ for all $x,y\in X$.
\end{lemma}

\begin{proof}
By induction on the construction.
The constant $1$ has values bounded by $1$ and Lipschitz constant $0$, and each $a_i$ has values bounded by $R$ and satisfies $|a_i(x)-a_i(y)|\leq d(x,y)$.
Rational multiples, sums, maxima and minima transform the two bounds exactly as they transform weights.
For a countable supremum, the elementary inequality $|\sup_nu_n-\sup_nv_n|\leq\sup_n|u_n-v_n|$ gives both bounds with $\sup_nw_n$, which is finite by the admission rule, and countable infima are treated in the same way.
\end{proof}

It is the weights, and not bounds on values, that keep the admitted suprema continuous.
A bound on the values of the terms alone would admit the supremum of the functions $x\mapsto\max(-1,\min(1,nx))$, for positive integers $n$, on the set $X=[-1,1]$ presented by $a_0(x)=x$ and $a_i=0$ for $i\geq1$; the supremum jumps from $0$ to $1$ at the origin (compare \cref{ex:overshoot}).
The function indexed by $n$ is admissible with weight $n$, and the weights have no common bound.

\begin{keytheorem}[Part~(1) of Theorem~A]\label{thm:an-representation}
Let $R$ be a positive real number, let $X$ be a nonempty set, and let $(a_i)_{i<\omega}$ be functions from $X$ to $[-R,R]$ that separate the points of $X$, with coordinate metric $d$.
If the metric space $(X,d)$ is separable, then the admissible functions are exactly the bounded Lipschitz functions on $(X,d)$.
\end{keytheorem}

\begin{proof}
By \cref{lem:an-weights}, every admissible function is bounded and Lipschitz.
Every real number $c$ in an interval $[-M,M]$, with $M$ a positive integer, is the supremum of a sequence of rational numbers $t_n\in[-M,M]$, so the constant function $c$ is admissible with weight at most $M$, as the supremum of the functions $t_n\cdot1$.
For an admissible function $g$ with weight $w$, the function $|g|=\max(g,-g)$ is admissible with weight $w$.

Let $f$ be bounded and Lipschitz, choose a positive integer $M>\norm f$, a positive rational number $L$ at least the Lipschitz constant of $f$, a positive integer $M_0\geq R$, and a dense sequence $(x_n)_{n<\omega}$ in $(X,d)$.
For each $n$, the function
\[
 d(\cdot,x_n)=\sup_{i<\omega}|a_i-a_i(x_n)|
\]
is admissible with weight at most $R+M_0$, and the functions $f(x_n)+L\,d(\cdot,x_n)$ are admissible with the common bound $M+L(R+M_0)$ on their weights.
Their infimum is therefore admissible, and the infimum equals $f$: the Lipschitz inequality gives $f(x)\leq f(x_n)+Ld(x,x_n)$ for every $n$, and for every $\epsilon>0$ an index $n$ with $d(x,x_n)<\epsilon$ gives $f(x_n)+Ld(x,x_n)\leq f(x)+2L\epsilon$.
\end{proof}

\begin{keytheorem}[Part~(2) of Theorem~A]\label{thm:an-approx}
Let $R$ be a positive real number, let $X$ be a nonempty set, and let $(a_i)_{i<\omega}$ be functions from $X$ to $[-R,R]$ that separate the points of $X$, with coordinate metric $d$.
If the metric space $(X,d)$ is separable, then the uniform closure of the set of admissible functions is $\BUC(X,d)$.
\end{keytheorem}

\begin{proof}
By \cref{lem:an-weights}, every uniform limit of admissible functions is bounded and uniformly continuous.
For the reverse inclusion, let $f\in\BUC(X,d)$ be given, and choose a positive integer $M>\norm f$.
For each positive integer $k$, define
\[
 g_k(x)=\inf_{y\in X}\bigl(f(y)+kd(x,y)\bigr).
\]
The functions inside the infimum are $k$-Lipschitz in $x$, and the infimum satisfies $-M\leq g_k\leq f\leq M$.
Hence the function $g_k$ is bounded and $k$-Lipschitz, and admissible by \cref{thm:an-representation}.

Let a positive number $\epsilon$ be given.
Uniform continuity gives a number $\delta>0$ such that $|f(x)-f(y)|<\epsilon$ whenever $d(x,y)<\delta$.
Choose $k$ with $k\delta\geq2M$.
If $d(x,y)<\delta$, then
\[
 f(y)+kd(x,y)\geq f(x)-\epsilon.
\]
If $d(x,y)\geq\delta$, then
\[
 f(y)+kd(x,y)\geq-M+k\delta\geq M>f(x).
\]
Taking the infimum over $y$ gives $f-\epsilon\leq g_k\leq f$ on all of $X$, for every sufficiently large $k$.
Thus $\norm{g_k-f}\xrightarrow{k\to\infty}0$.
\end{proof}

Every separable metric space carries a presentation to which Theorem~A applies.

\begin{corollary}[Separable metric spaces]\label{cor:polish}
Let $(X,\rho)$ be a nonempty separable metric space, for instance a Polish space, let $(x_i)_{i<\omega}$ be a dense sequence in $X$, and let $a_i(x)=\min(\rho(x,x_i),1)$ for $x\in X$ and $i<\omega$.
Then the functions $a_i$ take values in $[-1,1]$, separate the points of $X$, and have coordinate metric $\min(\rho,1)$.
Consequently, with $R=1$, the admissible functions are exactly the bounded Lipschitz functions on $(X,\rho)$, and the uniform closure of the set of these functions is $\BUC(X,\rho)$.
\end{corollary}

\begin{proof}
Write $d$ for the coordinate metric, and let $x,y\in X$ be given.
For every $i$, the numbers $a_i(x)$ and $a_i(y)$ lie in $[0,1]$.
The two numbers differ by at most $|\rho(x,x_i)-\rho(y,x_i)|\leq\rho(x,y)$, so $d(x,y)\leq\min(\rho(x,y),1)$.
Conversely, let $0<\epsilon<1$ be given, and choose $i$ with $\rho(x,x_i)<\epsilon$.
Then $a_i(x)<\epsilon$ and $a_i(y)\geq\min(\rho(x,y)-\epsilon,1)$, so $d(x,y)\geq\min(\rho(x,y),1)-2\epsilon$.
Hence $d=\min(\rho,1)$; in particular the functions $a_i$ separate points, and $(X,d)$ is separable.
A bounded function $f$ with Lipschitz constant $L$ for $\rho$ has Lipschitz constant at most $\max(L,2\norm f)$ for $\min(\rho,1)$, and a function that is Lipschitz for $\min(\rho,1)$ is Lipschitz for $\rho$; so the two metrics have the same bounded Lipschitz functions.
The two metrics agree on distances less than $1$, so they have the same bounded uniformly continuous functions.
The conclusion now follows from \cref{thm:an-representation,thm:an-approx}.
\end{proof}

\begin{example}[The rational numbers]\label{ex:rationals}
Completeness plays no role in Theorem~A: the map $x\mapsto(a_i(x))_{i<\omega}$ embeds $(X,d)$ isometrically into the cube $[-R,R]^\NN$ with the supremum metric, and the bounded Lipschitz and bounded uniformly continuous functions on $X$ are the restrictions of those on the completion of $X$, which is the closure of the image of $X$ in the cube.
For instance, present $\QQ$ as in \cref{cor:polish}, by an enumeration $(q_i)_{i<\omega}$ and $a_i(x)=\min(|x-q_i|,1)$.
The function $x\mapsto\min(|x-\sqrt2|,1)$, which measures the distance to a point missing from $\QQ$, is the admissible function $\inf_{i<\omega}(c_i+a_i)$, where $c_i=\min(|q_i-\sqrt2|,1)$.
The function that equals $1$ to the right of $\sqrt2$ and $0$ to the left of $\sqrt2$ is continuous on $\QQ$ but not uniformly continuous, and therefore lies outside $\BUC(\QQ,d)$: the theorem concerns the uniformity of $d$, not the topology of $X$.
\end{example}

On a compact space the countable operations are not needed, and the classical theorem already gives the conclusion.

\begin{remark}\label{rem:compact}
Suppose that the metric space $(X,d)$ is compact, and let $\mathcal A_0(X)$ be the set of admissible functions built without countable suprema and infima.
The set $\mathcal A_0(X)$ is a vector space over $\QQ$ that contains the rational constants and is closed under maxima and minima, and the members of $\mathcal A_0(X)$ are continuous by \cref{lem:an-weights}.
For distinct points $x$ and $y$ of $X$, real numbers $\alpha$ and $\beta$, and $\epsilon>0$, choose $i$ with $a_i(x)\ne a_i(y)$; rational numbers $s$ and $t$ with $|sa_i(x)+t-\alpha|<\epsilon$ and $|sa_i(y)+t-\beta|<\epsilon$ exist, and the function $sa_i+t$ belongs to $\mathcal A_0(X)$.
The lattice form of the Stone--Weierstrass theorem~\cite{Stone1948}, whose proof uses the two prescribed values only up to $\epsilon$, shows that $\mathcal A_0(X)$ is uniformly dense in the space of continuous functions on $(X,d)$, which on a compact space is $\BUC(X,d)$.
\end{remark}

On a particular space that is not compact the countable operations may still be superfluous, as on the interval $(0,1)$ presented by the identity function, but in general they cannot be dispensed with.

\begin{example}[A discrete space]\label{ex:discrete}
Let $R=1$ and $X=\NN$, and let $a_i(n)$ be $1$ if $n=i$ and $0$ otherwise.
The coordinate metric is the discrete metric, with distance~$1$ between distinct points, and the space $(\NN,d)$ is separable.
Every bounded function on $\NN$ is Lipschitz for this metric and is therefore admissible, by \cref{thm:an-representation}.
An admissible function built without countable suprema and infima, by contrast, uses only finitely many of the functions $a_i$, and is therefore constant outside a finite set.
Such a function differs from $n\mapsto(-1)^n$ by at least $1$ at some point, and the function $n\mapsto(-1)^n$ is therefore not a uniform limit of admissible functions built by finitely many operations.
\end{example}

\end{document}